\documentclass[a4paper]{article}

\usepackage{amsmath}
\usepackage{graphics}
\usepackage{latexsym}
\usepackage{amsfonts}
\usepackage{mathrsfs}
\usepackage{amssymb}
\usepackage{amssymb,amsmath,amsthm, amsfonts}

\numberwithin{equation}{section}

\newtheorem{propo}{Proposition}[section]

\newtheorem{lemma}{Lemma}[section]

\newtheorem{assumption}{Assumption}

\def\qed{ \ \vrule width.2cm height.2cm depth0cm\smallskip}

\def \ed {\end{document}}

\def \lb{\label}
\def\wr {w.r.t.}

\newcommand{\eps}{\varepsilon}

\newcommand{\brm}{\begin{rem}}
\newcommand{\ermq}{\end{rem}}

\newcommand{\ba}{\begin{array}}
\newcommand{\ea}{\end{array}}
\newcommand{\be}{\begin{equation}}
\newcommand{\ee}{\end{equation}}
\newcommand{\bea}{\begin{eqnarray}}
\newcommand{\eea}{\end{eqnarray}}
\newcommand{\beaa}{\begin{eqnarray*}}
\newcommand{\eeaa}{\end{eqnarray*}}

\def \R{I\!\!R}
\def \E{\mathbb{E}}

\def\a{\alpha}

\def\g{\gamma}
\def\d{\delta}

\def\z{\zeta}

\def\l{\lambda}

\def\si{\sigma}

\def\D{\Delta}

\def\cA{{\cal A}}
\def\cB{{\cal B}}
\def\cC{{\cal C}}

\def\cE{{\cal E}}
\def\cF{{\cal F}}

\def\cH{{\cal H}}

\def\cL{{\cal L}}

\def\cP{{\cal P}}

\def\cS{{\cal S}}

\def\no{\noindent}

\def\ms{\medskip}

\def\q{\quad}
\def\qq{\qquad}

\def \bp {\bold P}

\def\qed{ \hfill \vrule width.25cm height.25cm depth0cm\smallskip}

\newcommand{\basa}{\begin{assumption}}
\newcommand{\easa}{\end{assumption}}

\newcommand{\bas}{\begin{assum}}
\newcommand{\eas}{\end{assum}}

\def\dis{\displaystyle}

\def \P{\mathbb{P}}
\newtheorem{thm}{Theorem}[section]

\newtheorem{rem}[thm]{Remark}

\newtheorem{defn}[thm]{Definition}
\newtheorem{assum}[thm]{Assumption}

\newcommand{\rw}{\rightarrow}
\def \R{\mathbb{R}}

\def \esp {[0,T]\times \R^k}

\def \r {\mathbb{R}}
\def \rk {\r^k}
\def \mo {\{1,\dots,m\}}

\def \pg {\Pi_g}
\def \txp {(t,x)\in \esp}

\def \pgc {\Pi_g^c}
\title{Viscosity solutions for second order integro-differential equations without monotonicity condition: The probabilistic Approach} 
\author{
Said Hamad\`ene\thanks{Universit\'e du
Maine, LMM, Avenue Olivier Messiaen, 72085 Le Mans, Cedex 9, France, e-mail: hamadene@univ-lemans.fr} \,\, and \, Marie-Amelie Morlais\thanks{Universit\'e du Maine, LMM, Avenue Olivier Messiaen, 72085 Le Mans, Cedex 9, France, e-mail: Marie$_-$Amelie.Morlais@univ-lemans.fr} }

\begin{document}
\date{\today}
\maketitle

\begin{abstract}
In this paper, we establish a new existence and uniqueness result of a continuous viscosity solution for integro-partial differential equation (IPDE in short). 
 The novelty is that we relax the so-called monotonicity assumption on the driver which is classically assumed in the literature of viscosity solution of equation with non local terms. Our method strongly relies on the link between IPDEs and backward stochastic differential equations (BSDEs in short) with jumps for which we already know that the solution exists and is unique. In the second part of the paper, we deal with the IPDE with obstacle and we obtain similar results.
 \end{abstract}
\no{\bf AMS Classification subjects}: 
60H30
\medskip

\no {$\bf Keywords$}: Integro-differential equation ; Backward stochastic differential equation with jumps ; Viscosity solution ; Non-local term.
\section{Introduction}
In this paper, our objective is to establish a new existence and uniqueness result of the solution in viscosity sense of the following system of integro-partial differential equations: $\forall i\in \mo$, 
\be \left\{ \label{ipde-intro} \begin{array}{l} 
-\partial_t u^i(t,x)- b(t,x)^\top D_xu^i(t,x)-\frac{1}{2}\mathrm{Tr}\big(\sigma\sigma^\top (t,x) D^2_{xx}u^i(t,x)\big)-Ku^i(t,x)\\\\\q\qq -h^{(i)}(t,x, (u^i(t,x))_{i=1,m}, (\sigma^\top D_xu^i)(t,x), B_iu^i(t,x))=0,\,\, (t,x)\in \esp ;\\\\
u^i(T,x) = g^i(x),\; \;\forall \; i \in \{1,\cdots m\}, \; m \in \mathbb{N}^{*}
\end{array}  \right.
\ee where the operators $B_i$ and $K_i$ are defined by \be\label{defbk-intro}\begin{array}{l} B_iu^i(t,x)=\int_E\g^i(t,x,e)\left( u^i(t,x
+\beta(t,x,e))-u^i(t,x)\right)\l(de)\mbox{ and }\\\\
Ku^i(t,x)=\int_E\left( u^i(t,x
+\beta(t,x,e))-u^i(t,x)-\beta(t,x,e)^\top 
D_xu^i(t,x)\right)\l(de).\end{array}\ee  
We first note that, due to the presence of $B_iu^i$ and $K_iu^i$ in equation (\ref{ipde-intro}), such an IPDE is called of non-local type. IPDEs with non-local terms have been considered by several authors (see e.g. \cite{Alvareztourin,BarlesBuckPardoux,BarlesChassImb08,BarlesImbert08,sulem,HarrajOuknine05,Lundstrometal}, etc. and the references therein). It is by now well-known that this IPDE is connected with the following multi-dimensional backward stochastic differential equation with jumps: $\forall i\in \mo$, 
\be \label{bsde-jumps-intro}
   \left\{ \begin{array}{l} 
dY_s^{i;t,x} = -f^{(i)}(s,X^{t,x}_s,(Y^{i;t,x}_s)_{i=1,m},Z_s^{i;t,x}, U_s^{i;t,x})ds +Z_s^{i;t,x} dB_s +\int_{E} 
    U_{s}^{i;t,x}(e)\tilde{\mu}(ds,de),\,\,s\leq T\,;\\Y_T^{i;t,x}= g^{i}(X^{t,x}_T)\end{array}\right.\ee
where $\txp$, $B:=(B_s)_{s\leq T}$ is a Brownian motion, $\mu$ an independant Poisson random measure 
with compensator $ds\l(de)$ ($\l$ is the L\'evy measure of $\mu$) and $\tilde \mu (ds, de):=\mu (ds, de)-ds\l(de)$.
\ms 

For completeness, let us recall some already known results in the IPDE literature (and also those concerning the related BSDE with jumps). In \cite{TangLi94}, Tang-Li have shown that this BSDE with jumps (\ref{bsde-jumps-intro}) has a unique solution while Barles et al., in \cite{BarlesBuckPardoux}, have made the connection between this BSDE and the IPDE (\ref{ipde-intro}). Actually in \cite{BarlesBuckPardoux}, the authors have shown that if the coefficients $f^{(i)}$, $i=1,...,m$, have the following form:
\be\label{cond12}\begin{array}{l}
f^{(i)}(t,x,\vec{y},z,\z)= h^{(i)}(t,x, \vec{y}, z, \int_E\g_i(t,x,e)\z(e)\l(de))\end{array}
\ee
and, mainly, if  \\
\noindent 
(i) $\g_i\ge 0$ \\
(ii) $q\longmapsto h^{(i)}(t,x,\vec{y},z,q)$,  is non-decreasing ;\\
\no then the deterministic continuous functions $(u^i(t,x))_{i=1,m}$, defined by means of the representation of 
Feynman Kac's type of the processes $(Y^{i;t,x})_{i=1,m}$, i.e.,  
  \be \label{FKformula} \forall i=1,\dots,m,\,\,\displaystyle{ Y_{s}^{i;t,x} = u^{i}(s,X_{s}^{t,x}) \mbox{ for }s\in [t,T]\mbox{ and then }u^i(t,x):=Y^{i;t,x}_t},\ee
is the unique viscosity solution of (\ref{ipde-intro}) in the class of functions of polynomial growth. The two assertions (i)-(ii) above shall be referred later as the monotonicity conditions.
\ms

\no Therefore and in the first part of this paper, the main objective is to deal with IPDE (\ref{ipde-intro}) without assuming the two points (i)-(ii) above related to the non local term and the functions $h^{(i)}$. Actually we show that when the measure $\l$ is finite, equation (\ref{ipde-intro}) has a unique solution. Our method relies mainly on the following points:

(a) the characterization of the jump part of the BSDE (\ref{bsde-jumps-intro}) ;

(b) the existence and uniqueness of a solution of (\ref{bsde-jumps-intro}) for general $f^{(i)}$, $i=1,\dots,m$, which are merely Lipschitz in $(y,z,\zeta)$ and nothing more ; 

(c) the existence and uniqueness result of a solution of the IPDE (\ref{ipde-intro}) in the case
when $h^{(i)}$ does not depend on the component $\z$, which involves the jump part. This result is already obtained in 
\cite{BarlesBuckPardoux}. 
\ms

\noindent Thus, our main contribution consists in proving the following result:
there exists a unique viscosity solution of the general IPDE (\ref{ipde-intro}) (uniqueness
holds within the class of continuous functions with at most polynomial growth (w.r.t. $x$)). Besides,
the solution being given by the representation (\ref{FKformula}) of Feynman Kac's type, 
 we fill in the gap between existence and uniqueness results for BSDE with jumps of the form (\ref{bsde-jumps-intro}) (results which are already available for BSDEs and do not require the monotonicity conditions) and the results available
 in the IPDE literature. Finally, let us mention that, due to the presence of the operator $K$ inside the functional $h$ in IPDE (\ref{ipde-intro})
 and since neither (i) or (ii) holds, one cannot prove, as it is usual in viscosity literature,
 the classical comparison theorem. 
This motivates the introduction of a new definition of viscosity solution, which coincides  with the usual one when $h$ does not depend on the jump part and which allows to apply the classical results such as those obtained in \cite{BarlesBuckPardoux}. \\

\noindent According to the best of our knowledge and without assuming the two points (i)-(ii), 
such a result of existence and uniqueness of the solution of IPDE (\ref{ipde-intro}) has not been obtained so far.
We should also mention here one crucial point: since we assume the L\'evy measure $\l$ is bounded, the operators $B_iu^i$ are well-posed
 for functions which grow as polynomials $\wr$ $x$ at infinity. Thus we naturally introduce a new definition of viscosity solution (for IPDES with or without obstacle).
The main point is that even if our definitions are a bit
different from the ones given in \cite{BarlesBuckPardoux,sulem,HarrajOuknine05}, etc., but one can show that 
they coincide if (i)-(ii) above are satisfied. As a consequence, our study naturally extends the already known results in the IPDE literature.
\ms

\noindent In the second part of this paper, we consider the following IPDE with obstacle $(m=1)$:
\be \left\{ \label{obstacleipdeintro} \begin{array}{l} 
\min\Big \{u(t,x)-\ell (t,x); -\partial_t u(t,x)- b(t,x)^\top D_xu(t,x)-\frac{1}{2}\mathrm{Tr}\big(\sigma\sigma^\top (t,x) D^2_{xx}u(t,x)\big)\\\quad \qq-Ku(t,x)-h(t,x,u(t,x), (\sigma^\top D_xu)(t,x), Bu(t,x))\Big \}=0,\,\, (t,x)\in \esp ;\\
u(T,x) = g(x)
\end{array}  \right.
\ee
where the operators $Bu$ and $Ku$ are defined similarly as in (\ref{defbk-intro}) (just take $m=1$). Once again, this IPDE with obstacle (\ref{obstacleipdeintro}) is connected with the following reflected BSDE with jumps:
\begin{equation} \label{reflec-bsde-intro} \left\{ 
\begin{array}{l} 
dY^{t,x}_s = -f(s, X^{t,x}_s, Y^{t,x}_s,Z^{t,x}_s, U^{t,x}_s)ds  -dK^{t,x}_s+ Z^{t,x}_s dB_s +\int_E U^{t,x}_s(e) \tilde{\mu}(ds,de) , \,s\leq T;\\
  Y^{t,x}_s \ge \ell(s,X^{t,x}_s),\,\,s\le T \;\textrm{and} \; \int_{0}^{T}(Y^{t,x}_s -\ell (s,X^{t,x}_s) )dK^{t,x}_s =0 ;\\
 \displaystyle{Y^{t,x}_T =g(X^{t,x}_T)}\\
\end{array} \right.
 \end{equation}
for which Hamad\`ene-Ouknine \cite{hamouk} provide a unique solution for general generators $f$ by means of a fixed point theorem. The related IPDE is considered in several papers 
amongst one can quote (\cite{sulem, HarrajOuknine05}, etc.). However in those papers the conditions (i)-(ii) above on $\g_1$ and $h$ are assumed. Therefore our second main objective is to deal with the IPDE with obstacle (\ref{obstacleipdeintro}) for general functions $h$ and $\g$ which do not satisfy (i)-(ii). Indeed, similarly to the framework without obstacle, by using reflected BSDEs with jumps, we show that equation (\ref{obstacleipdeintro}) has a unique solution when the L\'evy measure $\l$ is bounded. This solution is also obtained with the help of the representation of Feynman Kac's formula of the unique solution of (\ref{reflec-bsde-intro}).
\ms

\noindent The outline of the paper is as follows: in the following second section, we provide all the 
necessary notations, assumptions and preliminary results concerning the study of general IPDEs (\ref{ipde-intro}) and related BSDEs with jumps as well.
In the third and fourth sections, we proceed with the two main results of the paper: (i) we first provide the main theoretical result of the paper,
i.e. the existence and uniqueness of the solution of the general non linear IPDE ; (ii) we generalize the result of the first part to IPDEs with obstacle.
 For completeness, usual definitions for viscosity solutions for both a non linear IPDE with and without obtacle are provided in an Appendix.
 
 \section{Preliminary results on BSDEs with jumps and their associated IPDEs}
\indent For sake of clarity, let us give the framework of our study as well as some notations which shall be used throughout the paper. In particular, we shall deeply rely on the relationship between the viscosity solution of some IPDEs and the solution of the related BSDE with jumps. Therefore and for sake of completeness, 
we need to introduce the stochastic framework and then give the connection with the integro-partial differential
equation we shall study.\\

\indent Let $(\Omega,\mathcal{F},(\mathcal{F}_t)_{t\leq T}, \mathbb{P})$ be a stochastic basis such that $\cF_0$ contains all $\P$-null sets of $\cF$, and $\cF_{t}=\cF_{t+}:=\bigcap_{\eps> 0}\cF_{t+\eps}$, $t\geq 0$, and we suppose that the filtration is generated by the two mutually independant processes:

\no (i) $B:=(B_t)_{t\ge 0}$ a $d$-dimensional Brownian motion and

\no (ii) a Poisson random measure $\mu$ on $\r^+\times E$, where $E:=\r^{\ell}-\{0\}$ is equipped with its Borel field $\cE$ ($\ell \ge 1$).
The compensator $\nu(dt,de)=dt\l(de)$ is such that $\{\tilde \mu([0,t]\times A)=(\mu-\nu)([0,t]\times A)\}_{t\geq 0}$ is a martingale for all $A\in \cE$  and 
satisfies $\l(A)<\infty$. We also assume that $\l$ is $\sigma$-finite measure on $(E,\cE)$ which integrates the function $(1\wedge |e|^2)_{e\in E}$.
\ms

Next we denote by: 

\no (iii) $\cP$ (resp. $\bp$) the field on $[0,T]\times \Omega$ 
of $(\cF_t)_{t\leq T}$-progressively measurable (resp. predictable) sets ;

\no (iv) $L^2(\l)$ the space of Borel measurable functions $\varphi:=(\varphi(e))_{e\in E}$ from $E$ into $\r$ such that $\|\varphi\|_{L^2(\l)}^2:=\int_E|\varphi(e)|^2\l(de)<\infty$ ; 

\no (v) $\mathcal{S}^{2}(\mathbb{R}^{\ell})$ $(\ell \in \mathbb{N}^*$) the space of RCLL (for right continuous with left limits) $\cP$-measurable and $\r^\ell$-valued processes such that $\mathbb{E}\big(\sup_{s\le T} |Y_s|^2\big) <\infty $ ; $\cA^2_c$ is its subspace of continuous non-decreasing processes $(K_t)_{t\leq T}$ such that $K_0=0$ ;  

\no (vi) $\mathcal{H}^2(\mathbb{R}^{\ell\times d})$ the space of processes $Z:=(Z_s)_{s\le T}$ which are $\cP$-measurable, $\r^{\ell \times d}$-valued and satisfying $\mathbb{E}[ \int_{0}^{T} |Z_s|^{2} ds ]<\infty $ ; 

\no (vii) 
$\mathcal{H}^2(L^2(\l))$ the space of processes $U:=(U_s)_{s\le T}$ which are $\bp$-measurable, $L^2(\l)$-valued and satisfying \\
$\mathbb{E}[ \int_{0}^{T} \|U_s(\omega)\|_{L^2(\l)}^{2} ds ]<\infty$ ; 

\no (viii) $\pg$ the set of deterministic functions $\varpi$: $(t,x)\in \esp \mapsto \varpi(t,x)\in \r$ of polynomial growth, i.e., for which there exists two constants $C$ and $p$ such that for any $(t,x)\in \esp$, 
$$
|\varpi(t,x)|\leq C(1+|x|^p).
$$The subspace of $\pg$ of continuous functions will be denoted by $\pgc$ ;

\no (ix) For any process $\theta:=(\theta_s)_{s\le T}$ and $t\in (0,T]$, $\theta_{t-}=\lim_{s\nearrow t}\theta_s$ and 
$\D_t \theta=\theta_t-\theta_{t-}$. 
\ms 

Now let $b$ and $\sigma$ be the following functions:
$$\ba{l}
b:(t,x)\in \esp \rw b(t,x) \in \r^k\\
\si :(t,x)\in\esp \rw \si (t,x)\in \r^{k\times d}.\ea
$$
We assume that they are jointly continuous in $(t,x)$ and Lipschitz continuous $\wr$ $x$ uniformly in $t$, i.e., there exists a constant $C$ such that  
\begin{equation}\label{bslip}
\forall \; (t,x, x') \in  [0,T] \times \mathbb{R}^{k+k},\;\; |b(t,x) -b(t,x'| +|\sigma(t,x) -\sigma(t, x')| \le C |x -x'|.
\end{equation}
Since $b$ and $\si$ are jointly continuous then by (\ref{bslip}), we easily deduce that they are of linear growth, i.e., there exists a constant $C$ such that  
\begin{equation}\label{bslip2}
\forall \; (t,x) \in [0,T] \times \mathbb{R}^d \;\;|b(t,x)| +|\sigma(t,x)|\le C(1+|x|).
\end{equation}

Let $\beta:(t,x,e)\in [0,T]\times \rk\times E \rw \beta (t,x,e)\in \rk$ be a measurable function such that for some real constant $C$, and for all $e\in E$, 
\begin{equation}\label{cdbeta}
\begin{array}{l}
\mbox{(i)} \q |\beta(t,x,e)|\le C (1\wedge |e|);\\
\mbox{(ii)}\q |\beta(t,x,e)-\beta (t,x',e)|\le C |x-x'|(1\wedge |e|);\\
\mbox{(iii)}\q \mbox{ the mapping }(t,x)\in \esp  \rw \beta (t,x,e)\in \rk \mbox{ is continuous uniformly }\wr \,\,e.
\end{array}
\end{equation}

Next let $(t,x)\in \esp$ and $(X^{t,x}_s)_{s\leq T}$ be the stochastic process solution of the following standard stochastic differential equation of diffusion-jump type:
\be \label{eq:diffusion}\begin{array}{l}X_{s}^{t,x} = x + \int_{t}^{s} b(r,X^{t,x}_r)dr +\int_{t}^{s}\sigma(r, X^{t,x}_r)dB_r
 +\int_{t}^{s}\int_{E}\beta(r, X^{t,x}_{r-}, e) \tilde{\mu}(dr,de), \mbox{ for } s\in [t,T] \mbox{ and }X_s^{t,x}=x \mbox{ if }s\leq t. \end{array}\ee
Under assumptions (\ref{bslip}), (\ref{bslip2}) and (\ref{cdbeta}) the solution of equation (\ref{eq:diffusion}) exists and is unique (see \cite{kunita} for more details). Moreover it satisfies the following estimates: $\forall p\geq 2$, $x,x'\in \r^k$ and $s\ge t$,
\be \label{estimx}
\E [\sup_{r\in [t,s]}|X^{t,x}_r-x|^p]\leq M_p(s-t)(1+|x|^p)] \mbox{ and }  
\E [\sup_{r\in [t,s]}|X^{t,x}_r-X^{t,x'}_r-(x-x')|^p]\leq M_p(s-t)|x-x'|^p
\ee
for some constant $M_p$.
\ms

Once for all and throughout this paper, we assume that (\ref{bslip}), (\ref{bslip2}) and (\ref{cdbeta}) are satisfied.
\ms

Let us now consider the following $m$-dimensional backward stochastic differential equation with jumps ($(t,x)\in \esp$): 
\be \left\{ \label{mainbsde}
    \begin{array}{l}  
\vec{Y}^{t,x}:=(Y^{i;t,x})_{i=1,m}\in \cS^2(\r^m),\,
{Z}^{t,x}\in \cH^2(\r^{m\times d}), U^{t,x}:=(U^{i;t,x})_{i=1,m}\in (\cH^2(L^2(\l))^m;\\ \forall i\in \mo, \,
Y_T^{i}= g^{i}(X^{t,x}_T) \mbox{ and }\forall s\leq T,\\
\qq dY_s^{i;t,x} = -f^{(i)}(s,X^{t,x}_s, \vec{Y}^{t,x}_s,Z_s^{i;t,x}, U_s^{i;t,x})ds -Z_s^{i;t,x} dB_s -\int_{E} 
    U_{s}^{i;t,x}(e)\tilde{\mu}(ds,de).\end{array}\right.  \ee
 where for any $i\in \mo$, 

(i) $f^{(i)}$ is a deterministic measurable function from $[0,T]\times \r^{k+m+m\times d}\times L^2(\l)$ into $\r$ ;            
 
 (ii) $Z_s^{i;t,x}$ is the $i$-th row of $Z_s^{t,x}$ and $U_s^{i;t,x}$ is the $i$-th component of $U_s^{t,x}$ ;
 
 (iii) $g^i$ are Borel measurable deterministic functions from $\r^k$ to $\r$. 
\ms

We now consider the following assumptions:
\ms

\noindent {\bf (H1)}: For any $i\in \mo$,
\ms

(i) $f^{(i)}$ is Lipschitz in $(y,z,u)$ uniformly in $(t,x)$, i.e., there exists a real constant $C$ such that for any $(t,x)\in \esp$, $(y,p,\zeta)$ and $(y',p',\zeta')$ elements of $\r^{m+d}\times L^2(\l)$, 
\be \label{cdlipschitzfyzu}   |f^{(i)}(t,x, y, p, \zeta) -f^{(i)}(t,x, y', p',\zeta')|
\leq C(|y-y'|+|p-p'|+\|\zeta-\zeta'\|_{L^2(\l)}). \ee 

(ii) the functions $f^{(i)}(t,x,0,0)$ and $g^i$ are of polynomial growth, i.e., belong to $\Pi_g$. 
\ms

\noindent {\bf (H2)}: For any $i\in \mo$:

(i) the functions $g_i$ are continuous ;
\ms

(ii) the mapping $(t,x)\in \esp \longmapsto f^{(i)}(t,x,\vec{y},z,\z)\in \r$ is continuous uniformly $\wr$ $(\vec{y},z,\z)$.

\bigskip

BSDEs with jumps have been already considered by Li-Tang in \cite{TangLi94} where they have provided the following result related to existence and uniqueness of the solution of (\ref{mainbsde}) (see also the paper by Barles et al. \cite{BarlesBuckPardoux}).
\ms

\begin{propo} \label{existencegene}(Tang-Li, \cite{TangLi94}): Assume that Assumption (H1) is fulfilled. Then for any $(t,x)\in \esp$, the BSDE (\ref{mainbsde}) has a unique solution 
$(\vec{Y}^{t,x},{Z}^{t,x},U^{t,x})$. \end{propo} 
            
Next let us consider the following structure condition on the functions $(f^{(i)})_{i=1,m}$.  
            \ms
            
            \no {\bf (H3)}: For any $i\in \mo$, there exists a Borel measurable deterministic function $h^{(i)}$ from 
            $[0,T]\times \r^{k+m+d+1}$ into $\r$ such that: 
             \be \label{structurecond}\begin{array}{l}f^{(i)}(t,x,\vec{y},z,\zeta)=h^{(i)}(t,x,\vec{y},z,\int_E\z (e)\g_i(t,x,e)\l(de)) \end{array}  \ee where for $i=1,\dots,m$, $\g_i$ is Borel measurable and verifies: $\forall (t,x,x')\in [0,T]\times \r^{k+k}$ and $e\in E$, there exists a constant $C\geq 0$, 
\be \label{condgamma}
\begin{array}{l}
(i)\q |\gamma_i(t,x,e)|\leq C(1\wedge |e|)\\
(ii)\q |\gamma_i(t,x,e)-\gamma_i(t,x',e)|\leq 
C(1\wedge |e|)|x-x'|
\\
(iii) \mbox{ the mapping } t\in [0,T]\longmapsto \g(t,x,e) \mbox{is continuous uniformly } \wr\,\, (x,e). \end{array}
\ee          

We then have the following result whose proof is given in Barles et al. (\cite{BarlesBuckPardoux}, Proposition 2.5 and 
Theorems 3.4, 3.5):
\begin{propo}\label{solvisco}(\cite{BarlesBuckPardoux}): Assume that (H1), (H2) and (H3) are fulfilled. Then there exist deterministic continuous functions $(u^i(t,x))_{i=1,m}$ which belong to $\pg$ such that for any $(t,x)\in \esp$, the solution of the BSDE (\ref{mainbsde}) verifies: 
\be\label{repre}
\forall i\in \mo,\,\,\forall s\in [t,T],\,\,Y^{i;t,x}_s=u^i(s,X^{t,x}_s).
\ee
(c) Moreover if for any $i\in \mo$, 

(c-i) $\g_i\geq 0$ ;

(c-ii) for any fixed $(t,x,\vec{y},z)\in [0,T]\times \r^{k+m+d}$, the mapping
$
q\in \r\longmapsto h^{(i)}(t,x,\vec{y},z,q)\in \r
$
is non-decreasing;\\ then $(u^i)_{i=1,m}$ is a continuous viscosity solution (in Barles et al.'s sense, see Definition \ref{bbpdef} in Appendix) of the following system of IPDEs: $\forall i\in \mo$, 
\be \left\{ \label{secondorder-pde} \begin{array}{l} 
-\partial_t u^i(t,x)- b(t,x)^\top D_xu^i(t,x)-\frac{1}{2}\mathrm{Tr}\big(\sigma\sigma^\top (t,x) D^2_{xx}u^i(t,x)\big)-Ku^i(t,x)\\\q\qq -h^{(i)}(t,x, (u^i(t,x))_{i=1,m}, (\sigma^\top D_xu^i)(t,x), B_iu^i(t,x))=0,\,\, (t,x)\in \esp ;\\
u^i(T,x) = g^i(x)
\end{array}  \right.
\ee where \be\label{defbk}\begin{array}{l} B_iu^i(t,x)=\int_E\g_i(t,x,e)\{u^i(t,x
+\beta(t,x,e))-u^i(t,x)\}\l(de)\mbox{ and }\\\\
Ku^i(t,x)=\int_E\{u^i(t,x
+\beta(t,x,e))-u^i(t,x)-\beta(t,x,e)^\top 
D_xu^i(t,x)\}\l(de).\end{array}\ee 
Finally, the solution $(u^i(t,x))_{i=1,m}$ is unique in the class of continuous functions of $\Pi_g$. 
\end{propo}
\begin{rem} By (\ref{repre}), for any $i\in \mo$ and $(t,x)\in \esp$,
\be\label{defui}
u^i(t,x):=Y^{i;t,x}_t.
\ee
\end{rem}
\section{The first main result: Existence and uniqueness of the solution for system of IPDEs}
To begin with, we are going to deal with the link between the stochastic process $U^{i;t,x}$ of the BSDE (\ref{mainbsde}) and the function $u^i$ defined in (\ref{defui}). For that, we need to assume additionally the following hypothesis on the L\'evy measure $\l$. 
\ms

\noindent {\bf (H4)}: The measure $\l$ is finite, i.e., $\l(E)<\infty$.
\ms 

\begin{propo} \label{caractu} 
Assume that (H1)-(H4) are fulfilled. Then for any $i=1,\dots,m$,
\be \label{caractu1}U_{s}^{i;t,x}(e)=
u^i(s,X^{t,x}_{s-}+\beta (s,X^{t,x}_{s-},e))-u^i(s,X^{t,x}_{s-}),\,\,ds\otimes d\mathbb{P}\otimes d\l \mbox{ on }[t,T]\times \Omega \times E.\ee
\end{propo}
\proof Let $i$ be fixed. First note that since $u^i$ belongs to $\pg$ and $\beta$ is bounded then by $\bf (H4)$ we have 
$$\E[\int_0^T\int_{E}\{|U^{i;t,x}_s(e)|^2+|u^i(s,X^{t,x}_{s-}+\beta(s,X^{t,x}_{s-},e)-
u^i(s,X^{t,x}_{s-})|^2\}\l(de)ds]<\infty$$ and hence, due to the finiteness of $\l$, one has
$$\E[\int_0^T\int_{E}\{|U^{i;t,x}_s(e)|+|u^i(s,X^{t,x}_{s-}+\beta(s,X^{t,x}_{s-},e)-
u^i(s,X^{t,x}_{s-})|\}\l(de)ds]<\infty.$$
Therefore and referring to \cite{rcont}, pp. 60, 
$$
\forall s\in [t,T],\,\,\int_t^{s}\int_{E} 
    U_{r}^{i;t,x}(e)\tilde{\mu}(dr,de)=\int_t^{s}\int_{E} 
    U_{r}^{i;t,x}(e){\mu}(dr,de)-\int_t^{s}\int_{E} 
    U_{r}^{i;t,x}(e)\l(de)dr.
$$On the other hand, since $Y^{i;t,x}$ satisfies the BSDE (\ref{mainbsde}) then for any $s\in [t,T]$,
\be \label{relation2}
\sum_{t<r\leq s}\{
Y^{i;t,x}_{r}-Y^{i;t,x}_{r-}\}=\int_t^{s}\int_{E} 
    U_{r}^{i;t,x}(e){\mu}(dr,de).
\ee Next, for any $s\in [t,T]$, $Y^{i;t,x}_s=u^i(s,X^{t,x}_s)$ and $u^i$ is continuous then 
\be \label{relation1}
\int_t^s\int_E(u^i(r,X^{t,x}_{r-}+\beta(s,X^{t,x}_{s-},e)-
u^i(r,X^{t,x}_{r-})\mu(dr,de)=\sum_{t<r\leq s}\{
Y^{i;t,x}_{r}-Y^{i;t,x}_{r-}\}.\ee It follows that for any $s\in [t,T]$,  
$$ \int_t^s\int_E(u^i(r,X^{t,x}_{r-}+\beta(r,X^{t,x}_{r-},e))-
u^i(r,X^{t,x}_{r-})-U_{r}^{i;t,x}(e))\mu(dr,de)=0.
    $$
Taking now the quadratic variation of this last process and then expectation to obtain $$
    \E[\int_t^Tdr\int_E|u^i(r,X^{t,x}_{r-}+\beta(r,X^{t,x}_{r-},e))-
u^i(r,X^{t,x}_{r-})-U_{r}^{i;t,x}(e)|^2\l(de)] =0
    $$which provides the desired equality.
\begin{rem}
This characterization of $U^{i;t,x}$ in terms of $u^i$ which is given in (\ref{caractu1}) plays a prominent role in the proof of our main result.
It is obtained under the condition (H4) of finiteness of the L\'evy measure $\l$. However
it can also be obtained under other conditions by using e.g. Malliavin calculus (see e.g. \cite{delong}, pp.84).
 But the use of Malliavin calculus requires stringent regularity condition on the data, therefore we do not use it as we are interested in obtaining results for quite general IPDEs. 
 \end{rem}
\subsection{Existence and uniqueness of the solution of the system of IPDEs}
We first give our meaning of the definition of the viscosity solution of system (\ref{ipde-intro}). It is not exactly the same as the one of Barles et al.'s paper (see Definition \ref{bbpdef} in Appendix). 

For any function $\phi$ belonging to $\cC^{1,2}([0,T]\times \r^k)$ and $\r$-valued, we define $\cL^X\phi$ by
$$
\cL^X\phi(t,x)=\frac{1}{2}\mbox{Tr}[(\sigma\sigma^\top)(t,x)D^2_{xx}\phi(t,x)]+b(t,x)^\top D_x\phi(t,x)
+K\phi(t,x), \, (t,x)\in \esp,
$$
where $K\phi(t,x)$ is given in (\ref{defbk}), and it is actually well-posed for any $\phi$ in $\cC^{1,2}([0,T]\times \r^k)$.

\begin{defn}\label{od}
A family of deterministic functions $u =(u^i)_{i=1,m}$, such that, for any $i\in \mo$, the map $u^i: (t,x) \mapsto u^i(t,x)$ belongs to $\Pi_{g}^c$ (spaces of continuous functions with at most polynomial growth w.r.t $x$), is said to be a viscosity sub-solution (resp. super-solution)
 of the IPDE (\ref{ipde-intro}) if: $\forall i\in \mo$, \\
(i) $\forall x\in \r^k$, $u_i(T, x) \le g_i(x) $ (resp. $u_i(T, x) \ge g_i(x) $) ;\\
 (ii)  For any  $(t,x)\in (0,T)\times \r^k$ and any function $\phi$ of class $\cC^{1,2}(\esp)$ such that $(t,x)$ is a global maximum (resp. minimum) point of 
$ u_i -\phi$ and $(u_i-\phi)(t,x)=0$, one has
$$ -\partial_t \phi(t,x)
-\mathcal{L}^X \phi(t,x) -h^{(i)}(t,x, (u^i(t,x))_{i=1,m}, \sigma^\top(t,x)D_x\phi(t,x), B_i(u^{i})(t,x)) \le 0, 
$$
(resp. 
$$-\partial_t \phi(t,x)
-\mathcal{{L}}^X \phi(t,x) -h^{(i)}(t,x, (u^i(t,x))_{i=1,m}, \sigma^\top(t,x)D_x\phi(t,x), B_i(u^{i})(t,x)) \ge 0.)
$$

The family $u =(u^i)_{i=1,m}$ is a viscosity solution of (\ref{ipde-intro}) if it is both a viscosity sub-solution and viscosity super-solution. 
\end{defn}

Let us now compare the two definitions \ref{od} and \ref{bbpdef} of viscosity solutions:

\begin{rem}\label{coincid} ${}$

\no (i) If for any $i\in \mo$, the function $h^{(i)}$ does not depend on its last component $\z$ then Definitions 
\ref{od} and \ref{bbpdef} are the same. 
\ms

\no (ii) In our Definition \ref{od}, we have used $B_iu^i$ instead of $ B_{i} \phi$: indeed $B_i u^{i}$
is well posed since $u^i$ is in $ \pg$, $\beta$ is bounded and $\l$ finite while it is replaced by $B_i\phi$ ($\phi$ being the 
\textit{smooth} test function) in Barles et al.'s definition (Definition \ref{bbpdef} in Appendix).
In the latter definition of Barles et al.'s this is the lack of regularity of $u^i$ which makes that $B_iu^i$ is ill-posed. 

\noindent (iii) We finally mention that this new definition is crucial in our proof of the existence result and more precisely to prove that our candidate is a viscosity solution:
 more precisely and by using this new definition,
 we shall be able to use classical existence results for IPDE obtained in Barles et al. \cite{BarlesBuckPardoux}: the main point being that,
 when the non local term is frozen, we fall in Barles et al.'s framework. 
\end{rem}

The main result of this paper is the following one:
\begin{thm} \lb{mainthm} Under Assumptions (H1)-(H4), the family of functions $(u^i)_{i=1,m}$ which belong to $\pgc$
defined in (\ref{defui}) is a viscosity solution of (\ref{ipde-intro}). Moreover it is unique in the class $\pgc$. 
\end{thm}
\proof Let us consider the following multi-dimensional BSDE:
\be \left\{ \label{mainbsde2}
    \begin{array}{l}  
\vec{\underbar Y}^{t,x}:=(\underbar Y^{i;t,x})_{i=1,m}\in \cS^2(\r^m),\,
{\underbar Z}^{t,x}\in \cH^2(\r^{m\times d}),\, \underbar U^{t,x}:=(\underbar U^{i;t,x})_{i=1,m}\in (\cH^2(L^2(\l))^m;\\\\ \forall i\in \mo,\,\,
\underbar Y_T^{i}= g^{i}(X^{t,x}_T) \mbox{ and }\forall s\le T, \\\\
\qq d\underbar Y_s^{i;t,x} = -h^{(i)}(s,X^{t,x}_s, \vec{\underbar Y}^{t,x}_s,\underbar Z_s^{i;t,x}, \int_E
\g_i(s,X^{t,x}_s,e)\{u^i(s,X^{t,x}_s+\beta (s,X^{t,x}_s,e))-u^i(s,X^{t,x}_s)\}\l(de))ds\\\\ \qq\qq\qq +\underbar Z_s^{i;t,x} dB_s +\int_{E} 
    \underbar U_{s}^{i;t,x}(e)\tilde{\mu}(ds,de).\end{array}\right.  \ee
As, for any $i=1,...,m$, $u^i$  belongs to $\pgc$, $\beta(t,x,e)$ is bounded and verifies (\ref{cdbeta}) and since $\l$ is finite,
the solution of this equation exists and is unique by Proposition \ref{existencegene}, noting that the functions $g^i$ and 
$$
(t,x,y,z)\longmapsto h^{(i)}(t,x,y,z,
\int_E
\g_i(t,x,e)\{u^i(t,x+\beta (t,x,e))-u^i(t,x)\}\l(de))$$ verify Assumptions (H1). Moreover, by Proposition \ref{solvisco}, there exists a family of deterministic continuous functions of polynomial growth $(\underbar u^i)_{i=1,m}$ such that for any $(t,x)\in \esp$,
$$
\forall s\in [t,T],\,\, \underbar Y^{i;t,x}_s=
\underbar u^i(s,X^{t,x}_s).
$$
Next by Proposition \ref{solvisco} and Remark \ref{coincid}-(i), the family $(\underbar u^i)_{i=1,m}$ is a viscosity solution of the following system:
\be \left\{ \label{secondorder-pde-inter} \begin{array}{l} 
-\partial_t \underbar u^i(t,x)- b(t,x)^\top D_x\underbar u^i(t,x)-\frac{1}{2}\mathrm{Tr}\big(\sigma\sigma^\top (t,x) D^2_{xx}\underbar u^i(t,x)\big)-K\underbar u^i(t,x)\\\q\qq -h^{(i)}(t,x, (\underbar u^i(t,x))_{i=1,m}, (\sigma^\top D_x\underbar u^i)(t,x), B_iu^i(t,x))=0,\,\, (t,x)\in \esp ;\\
\underbar u^i(T,x) = g^i(x)
\end{array}  \right.
\ee
Note that in this system (\ref{secondorder-pde-inter}), the last component of $h^{(i)}$ is $B_iu^i(t,x)$ and not $B_i\underbar u^i(t,x)$.

Next and once more, let us consider the system of BSDEs by which the family $(u^i)_{i=1,m}$ is defined through the Feynman Kac's formula (\ref{defui}).
Such a system of BSDEs is given by 
\be \left\{ \label{mainbsde21}
    \begin{array}{l}  
\vec{ Y}^{t,x}:=( Y^{i;t,x})_{i=1,m}\in \cS^2(\r^m),\,
{ Z}^{t,x}\in \cH^2(\r^{m\times d}),\,  U^{t,x}:=( U^{i;t,x})_{i=1,m}\in (\cH^2(L^2(\l))^m;\\\\ \forall i\in \mo,\,\,
 Y_T^{i}= g^{i}(X^{t,x}_T) \mbox{ and }\forall s\le T,\\\\
\qq dY_s^{i;t,x} = -h^{(i)}(s,X^{t,x}_s, \vec{Y}^{t,x}_s,Z_s^{i;t,x}, \int_E
\g_i(s,X^{t,x}_s,e) U_{s}^{i;t,x}(e)\l(de))ds\\\\ \qq\qq\qq +Z_s^{i;t,x} dB_s +\int_{E} 
    U_{s}^{i;t,x}(e)\tilde{\mu}(ds,de).\end{array}\right.  \ee
Since we know that, for any $i$ in $\{1, \cdots,m\}$, $u^i$ belongs to $\pgc$, therefore and due to Proposition \ref{caractu}, one has
$$
U_{s}^{i;t,x}(e)=
u^i(s,X^{t,x}_s+\beta (s,X^{t,x}_s,e))-u^i(s,X^{t,x}_s),\,\,ds\otimes d\mathbb{P}\otimes d\l \mbox{ on }[t,T]\times \Omega \times E.
$$
Plugging now this relation in the first term of the right-hand side of the second equality of (\ref{mainbsde21}), one obtains, by uniqueness of the solution of the BSDE (\ref{mainbsde2}), that for any $s\in [t,T] \mbox{ and }i\in \mo,\,\,
\underbar Y^{i;t,x}_s=Y^{i;t,x}_s$. Thus for any $i\in \mo$, $u^i=\underbar u^i$. Henceforth, the family $(u^i)_{i=1,m}$ is a viscosity solution of (\ref{ipde-intro}) in the sense of Definition \ref{od}. \\

Next, let us show that it is unique in the class $\pgc$. So let $(\bar u^i)_{i=1,m}$ be another family of $\pgc$ which is solution of the system (\ref{ipde-intro}) in the sense of Definition \ref{od}. 
 
Let us consider the following system of BSDEs:
\be \left\{ \label{mainbsde3}
    \begin{array}{l}  
\vec{\bar Y}^{t,x}:=(\bar Y^{i;t,x})_{i=1,m}\in \cS^2(\r^m),\,
{\bar Z}^{t,x}\in \cH^2(\r^{m\times d}),\, \bar U^{t,x}:=(\bar U^{i;t,x})_{i=1,m}\in (\cH^2(L^2(\l))^m;\\\\ \forall i\in \mo,\,\,
\bar Y_T^{i}= g^{i}(X^{t,x}_T) \mbox{ and }\forall s\le T,\\\\
\qq d\bar Y_s^{i;t,x} = -h^{(i)}(s,X^{t,x}_s, \vec{\bar Y}^{t,x}_s,\bar Z_s^{i;t,x}, \int_E
\g_i(s,X^{t,x}_s,e)\{\bar u^i(s,X^{t,x}_s+\beta (s,X^{t,x}_s,e))-\bar u^i(s,X^{t,x}_s)\}\l(de))ds\\\\ \qq\qq\qq +\bar Z_s^{i;t,x} dB_s +\int_{E} 
    \bar U_{s}^{i;t,x}(e)\tilde{\mu}(ds,de).\end{array}\right.  \ee
Therefore there exists a family of deterministic continuous functions $(v^i)_{i=1,m}$ of class $\Pi_g$ such that 
$$
\forall s\in [t,T],\,\, \bar Y^{i;t,x}_s=v^i(s,X^{t,x}_s).
$$
Additionally, by Definition \ref{bbpdef} and Proposition \ref{solvisco}, $(v^i)_{i=1,m}$ is the unique solution in the subclass $\pgc$ of continuous functions of the following system:

\be \left\{ \label{secondorder-pde-inter2} \begin{array}{l} 
-\partial_t v^i(t,x)- b(t,x)^\top D_xv^i(t,x)-\frac{1}{2}\mathrm{Tr}\big(\sigma\sigma^\top (t,x) D^2_{xx}v^i(t,x)\big)-Kv^i(t,x)\\\q\qq -h^{(i)}(t,x, (v^i(t,x))_{i=1,m}, (\sigma^\top D_xv^i)(t,x), B_i\bar u^i(t,x))=0,\,\, (t,x)\in \esp ;\\
v^i(T,x) = g^i(x)
\end{array}  \right.
\ee
But, the family $(\bar u^i)_{i=1,m}$ belongs to $\pgc$ and solves system (\ref{secondorder-pde-inter2}). Therefore, by the uniqueness result of Proposition \ref{solvisco} and Remark \ref{coincid} -(i), for any $i\in \mo$, one deduces that $\bar u^i= v^i$. On the other hand, by the 
characterization of the jumps of Proposition \ref{caractu}, for any 
$i\in \mo$, it holds
$$\begin{array}{ll}
\bar U_{s}^{i;t,x}(e)&=
v^i(s,X^{t,x}_{s-}+\beta (s,X^{t,x}_{s-},e))-v^i(s,X^{t,x}_{s-})\\{}&=\bar u^i(s,X^{t,x}_{s-}+\beta (s,X^{t,x}_{s-},e))-\bar u^i(s,X^{t,x}_{s-}).\ea$$Next by replacing in (\ref{mainbsde3}) the quantity $\bar u^i(s,X^{t,x}_{s-}+\beta (s,X^{t,x}_{s-},e))-\bar u^i(s,X^{t,x}_{s-})$ with $\bar U_{s}^{i;t,x}(e)$, we deduce that 
$$ \forall \; i\in \mo, \quad \quad
\bar Y^{i;t,x}=Y^{i;t,x}$$  since the solution of the BSDE (\ref{mainbsde3}) is unique. Thus and for any $i\in \mo$, $u^i=\bar u^i=v^i$ which means that the solution of (\ref{ipde-intro}) in the
sense of Definition \ref{od} is unique inside the class $\pgc$. 
\begin{rem} Since the L\'evy measure $\l$ is assumed to be bounded then one can
relax a bit the conditions (\ref{cdbeta}) and (\ref{condgamma}) on $\beta$ and $(\gamma_i)_{i=1,m}$ respectively. 
\end{rem}
\subsection{The second main result: Generalization to IPDEs with obstacles }\label{finalremark}
The previous result can be generalized to IPDEs with one (either lower or upper) obstacle. Actually assume that $m=1$ and let us denote $f^{(1)}$, $h^{(1)}$, $g^1$ and $\g_1$ simply by $f$, $h$, $g$ and $\g$ respectively. Next let us consider the following IPDE with obctacle $\ell$, which is a function of $(t,x)$:
\be \left\{ \label{obstacleipde} \begin{array}{l} 
\min\Big \{u(t,x)-\ell (t,x); -\partial_t u(t,x)- b(t,x)^\top D_xu(t,x)-\frac{1}{2}\mathrm{Tr}\big(\sigma\sigma^\top (t,x) D^2_{xx}u(t,x)\big)\\\quad \qq-Ku(t,x)-h(t,x,u(t,x), (\sigma^\top D_xu)(t,x), Bu(t,x))\Big \}=0,\,\, (t,x)\in \esp ;\\
u(T,x) = g(x)
\end{array}  \right.
\ee
where once again the operators $Bu$ and $Ku$ are given by:
\be\label{defbkr}\begin{array}{l} Bu(t,x)=\int_E\g(t,x,e)\{u(t,x
+\beta(t,x,e))-u(t,x)\}\l(de)\mbox{ and }\\\\
Ku(t,x)=\int_E\{u(t,x
+\beta(t,x,e))-u(t,x)-\beta(t,x,e)^\top 
D_xu(t,x)\}\l(de).\end{array}\ee 
Note that under (H4) if $u$ belongs to $\pg$ then the operator $Bu$ is well-posed.
\ms

The general reflected BSDE with jumps associated with IPDE with obstacle (\ref{obstacleipde}) is the following one: 
\begin{equation} \label{reflec-bsde} \left\{ 
\begin{array}{l} 
{Y}^{t,x}\in \cS^2(\r),\,
{Z}^{t,x}\in \cH^2(\r^d), U^{t,x}\in \cH^2(L^2(\l))  \mbox{ and }K^{t,x}\in \cA^2_c\,;\\\\
 dY^{t,x}_s = -f(s, X^{t,x}_s, Y^{t,x}_s,Z^{t,x}_s, U^{t,x}_s)ds  -dK^{t,x}_s+ Z^{t,x}_s dB_s +\int_E U^{t,x}_s(e) \tilde{\mu}(ds,de) , \,s\leq T;\\\\
  Y^{t,x}_s \ge \ell(s,X^{t,x}_s),\,\,s\le T \;\textrm{and} \; \int_{0}^{T}(Y^{t,x}_s -\ell (s,X^{t,x}_s) )dK^{t,x}_s =0 ;\\\\
 \displaystyle{Y^{t,x}_T =g(X^{t,x}_T)}\\\\
\end{array} \right.
 \end{equation}
where $(t,x)\in \esp$ is fixed. 

The following result related to existence and uniqueness of a solution for this reflected BSDE with jumps (\ref{reflec-bsde}) is given in (\cite{hamouk}, Theorem 1.2.b).
\begin{propo} \label{propinter1}\cite{hamouk} Assume that:

\no (i) $f$ is Lipschitz in $(y,z,\z)\in \r^{1+d}\times L^2(\l)$ uniformly w.r.t. $(t,x)$ and the function \\ $(t,x)\in \esp \longmapsto f(t,x,0,0,0)$ belongs to $\Pi_g$ ; 

\no (ii) $g$ belongs to $\Pi_g$ and $\ell (T,x)\geq g(x)$, $\forall x\in \r$ ; 

\no (iii) $\ell$ is continuous and belongs to $\Pi_g$.  
\medskip

\no Then the BSDE (\ref{reflec-bsde}) has a unique solution $(Y^{t,x},Z^{t,x},U^{t,x},K^{t,x})$. Moreover is satisfies the following estimate:
\be\label{estimyzv}
\E\Big [\sup_{s\leq T}|Y^{t,x}_s|^2+(K^{t,x}_T)^2+\int_0^T
\{|Z^{t,x}_s|^2+\|U^{t,x}_s\|_{L^2(\l)}^2\}ds\Big ]\le C
\E\Big [|g(X^{t,x}_T)|^2+
\sup_{s\leq T}|\ell(s,X^{t,x}_s)|^2+\int_0^T|f(s,X^{t,x}_s,0,0,0)|^2ds\Big ].
\ee
\end{propo}

Our main objective is now to connect the solution of the RBSDE with jumps with the solution in viscosity sense of IPDE with obstacle (\ref{obstacleipde}). To begin with let us precise the definition of viscosity solution we deal with.

\begin{defn}\label{od2}
We say that a function $u(t,x)$ which belongs to $\pgc$ is a viscosity sub-solution (resp. super-solution) of the IPDE (\ref{obstacleipde}) if: \\
(i) $\forall x\in \r^k$, $u(T, x) \le g(x) $ (resp. $u(T, x) \ge g(x) $) ;\\
 (ii)  For any  $(t,x)\in (0,T)\times \r^k$ and any function $\phi$ of class $\cC^{1,2}(\esp)$ such that $(t,x)$ is a global maximum (resp. minimum) point of 
$ u -\phi$ and $(u-\phi)(t,x)=0$, one has
$$\min\Big\{u(t,x)-\ell(t,x); -\partial_t \phi(t,x)
-\mathcal{L}^X \phi(t,x) -h(t,x, u(t,x), \sigma^\top(t,x)D_x\phi(t,x), Bu(t,x)) \Big \}\le 0, 
$$
(resp. 
$$\min\Big\{u(t,x)-\ell(t,x);-\partial_t \phi(t,x)
-\mathcal{{L}}^X \phi(t,x) -h(t,x, u(t,x), \sigma^\top(t,x)D_x\phi(t,x), Bu(t,x))\Big\} \ge 0.)
$$
The function $u$ is a viscosity solution of (\ref{obstacleipde}) if it is both a viscosity sub-solution and viscosity super-solution. 
\end{defn}
Next let us introduce the following assumptions:
\ms

\no {\bf (H5)}

\no (i) The assumptions of Proposition \ref{propinter1} are satisfied ;
 
\no (ii) the function $\g(t,x,e)$ verifies 
(\ref{condgamma}) ;

\no (iii) the function $h(t,x,y,z,\eta)$ such that 
\be \label{structurecond}\begin{array}{l}f(t,x,y,z,\zeta)=h(t,x,y,z,\int_E\z (e)\g(t,x,e)\l(de)) \end{array}  \ee
is continuous in $(t,x,y,z,\eta)$ and Lipschitz in $(x,y,z,\eta)$ uniformly w.r.t. $t$ ;

\no (iv) the function $g$ is continuous in $x$. 
\ms

We then have the following result related to the solution of (\ref{reflec-bsde}) which exists under (H5). 
\ms

\begin{propo} Assume that (H4)-(H5) are fulfilled. Then there exists a continuous deteministic function $u$ which belongs to $\pgc$ such that:
\be\label{fctu}
\forall s\in [t,T], Y^{t,x}_s=u(s,X^{t,x}_s).
\ee and 
\be \label{caracsaut}U_{s}^{t,x}(e)=
u(s,X^{t,x}_{s-}+\beta (s,X^{t,x}_{s-},e))-u(s,X^{t,x}_{s-}),\,\,ds\otimes d\mathbb{P}\otimes d\l \mbox{ on }[t,T]\times \Omega \times E.\ee
\end{propo}
\proof
Let $\Sigma:=\cH^2(\r)\times \cH^2(L^2(\l))$ and $\Psi$ be the functional which with a pair of processes $(y,v)$ which belongs to $\Sigma$ associates 
$\Psi(y,v):=(Y,V)$ such that $(Y,Z,V,K)$ is the solution of the following reflected BSDE with jumps: 
\begin{equation} \label{reflec-bsde2} \left\{ 
\begin{array}{l} 
{Y}\in \cS^2(\r),\,
{Z}\in \cH^2(\r^d), V\in \cH^2(L^2(\l))
\mbox{ and } K\in \cA_c^2\,\,;\\\\
 dY_s = -f(s, X^{t,x}_s, y_s,Z_s, v_s)ds  -dK_s+ Z_s dB_s +\int_E V_s(e) \tilde{\mu}(ds,de) , \,s\leq T;\\\\
  Y_s \ge \ell(s,X^{t,x}_s),\,\,s\le T \;\textrm{and} \; \int_{0}^{T}(Y_s -\ell (s,X^{t,x}_s) )dK_s =0 \,;\\\\
 \displaystyle{Y_T =g(X^{t,x}_T)}\\\\
\end{array} \right.
 \end{equation}
where $(t,x)\in \esp$ is fixed (we have omitted the dependance in $(t,x)$ of $(Y,V)$ as there is no confusion). The solution of this equation exists thanks to (\cite{hamouk}, Theorem 1.2.b). Next for $\a\in \r$, let us define the norm $\|.\|_\a$ on $\Sigma$ by:
$$
\|(y,v)\|_\a:=\sqrt{\E[\int_0^Te^{\a s}\{|y_s|^2+\|v_s\|_{L^2(\l)}^2\}ds]}.$$
As in \cite{hamouk}, Theorem 1.2.b, one can show that for an appropriate $\a_0$, $\Psi$ is a contraction on $(\Sigma,\|.\|_{\a_0})$ and thus, this mapping has a unique fixed point $(Y^{t,x},U^{t,x})$ which with $Z^{t,x}$ and $K^{t,x}$ gives the unique solution of (\ref{reflec-bsde}). Let us now consider the following sequence of processes:
$$
(Y^0,V^0)=(0,0) \mbox{ and for }n\geq 1,\,\,
(Y^{n},V^{n})=
\Psi(Y^{n-1},V^{n-1}).
$$ 
Then obviously $(Y^{n},V^{n})_{n\ge 0}$ converges in $(\Sigma,\|.\|_{\a_0})$ to 
$(Y^{t,x},U^{t,x})$. Next by an induction argument we have: for any $n\geq 0$:

(i) there exists a deterministic continuous function $u^n:(t,x)\in \esp \mapsto u^n(t,x)$ which belongs to $\Pi_g$ such that for any $s\in [t,T]$, $Y^n_s=u^n(s,X^{t,x}_s)$ ;

(ii) $u^n$ belongs to $\pgc$ and 
$$
V^{n}_s(e):=
u^{n}(s,X^{t,x}_{s-}+\beta (s,X^{t,x}_{s-},e))-u^{n}(s,X^{t,x}_{s-}),\,\,ds\otimes d\mathbb{P}\otimes d\l \mbox{ on }[t,T]\times \Omega \times E.
$$

\no Indeed for $n=0$, the properties (i), (ii) are valid. So suppose that they are satified for some $n$. Then $(Y^{n+1},Z^{n+1},V^{n+1},K^{n+1})$ verifies:
$\forall s\in [t,T]$, 
\begin{equation} \label{reflec-bsde3} \left\{ 
\begin{array}{l} 
 dY^{n+1}_s = -f(s, X^{t,x}_s, u^n(s,X^{t,x}_s),Z^{n+1}_s,\int_E\{
u^{n}(s,X^{t,x}_{s-}+\beta (s,X^{t,x}_{s-},e))-u^{n}(s,X^{t,x}_{s-})\}\g (s,X^{t,x}_{s-},e)\l (de))ds  \\\\\qq\qq\qq-dK^{n+1}_s+ Z^{n+1}_s dB_s +\int_E V^{n+1}_s(e) \tilde{\mu}(ds,de),\,\,;\\\\
  Y^{n+1}_s \ge \ell(s,X^{t,x}_s)\textrm{ and } \; \int_{t}^{T}(Y^{n+1}_s -\ell (s,X^{t,x}_s) )dK^{n+1}_s =0 \,;\\\\
 \displaystyle{Y^{n+1}_T =g(X^{t,x}_T).}
\end{array} \right.
 \end{equation}
Therefore the existence and continuity of $u^{n+1}$ are obtained in the same way as in 
(\cite{sulem}, \cite{HarrajOuknine05}) since the generator of $Y^{n+1}$ does not depend on $V^{n+1}$. Note that by (\ref{estimyzv}) we easily deduce that $u^{n+1}$ belongs to $\pg$. 
Finally the last property of (ii), i.e.,
$$
V^{n+1}_s(e):=
u^{n+1}(s,X^{t,x}_{s-}+\beta (s,X^{t,x}_{s-},e))-u^{n+1}(s,X^{t,x}_{s-}),\,\,ds\otimes d\mathbb{P}\otimes d\l \mbox{ on }[t,T]\times \Omega \times E.
$$
is obtained in a similar fashion as in Proposition \ref{caractu}. The proof of the induction procedure and thus of the two claims (i) and (ii) is now complete.
\ms

It now remains to justify that the two representations (\ref{fctu}) and (\ref{caracsaut}) of both processes $Y^{t,x}$ and $U^{t,x}$  hold at the limit (when $n$ goes to $\infty$). Then and to proceed we note that the following inequality holds:
\be\lb{inegdesk}\ba{l}
\int_t^T( Y^{n+1;t,x}_s-Y^{n+1;t,x}_s)d(
K^{n+1;t,x}_s-K^{n+1;t,x}_s)\leq 0.\ea \ee Next by It\^o's formula and (\ref{inegdesk}) for any $n,m\geq 0$ we have: $\forall s\in [t,T]$, 
$$\ba{l}
(Y^{n+1}_s-Y^{m+1}_s)^2+\int_s^T
|Z^{n+1}_r-Z^{m+1}_r|^2dr+\sum_{s<r\leq T}(\Delta_r(Y^{n+1}-Y^{m+1}_r))^2\\\\\q \leq \int_s^T(Y^{n+1}_r-Y^{m+1}_r)(
f(r, X^{t,x}_r,Y^{n}_r ,Z^{n+1}_r,V^n_r)-
f(r, X^{t,x}_r,Y^{m}_r ,Z^{m+1}_r,V^m_r))dr\\\\\q 
-2\int_s^T(Y^{n+1}_r-Y^{m+1}_r)(Z^{n+1}_r-Z^{m+1}_r)dB_r
-2\int_s^T\int_E\{(Y^{n+1}_{r-}-
Y^{m+1}_{r-})(V^{n+1}_r(e)-V^{m+1}_r(e))\tilde \mu(de,dr).\ea
$$
Then in a classical way we obtain that 
$$
\E[\sup_{t\le s\le T}|Y^{n+1}_s-Y^{m+1}_s|^2]\rw 0
\mbox{ as }n,m\rw \infty.
$$Therefore, the sequence of functions $(u^n)_{n\ge 0}$ converges pointwisely in $[0,T]\times \r^k$ to a deterministic function $u$. Moreover for any $(t,x)\in \esp$ we have 
$$\forall s\in [t,T],\,
Y^{t,x}_s=u(s,X^{t,x}_s).$$
Finally the continuity of $u$ is obtained in a similar way as in ((\cite{sulem}), pp.6) or (\cite{HarrajOuknine05}), pp.45) and relying to the proof of Proposition \ref{caractu} we obtain 
(\ref{caracsaut}) since $\l$ is finite. \qed
\ms

We now are ready to give the main result of this subsection.
\begin{thm}Assume that (H4), (H5) are fulfilled. Then the function $u$ defined in (\ref{fctu}) is the unique viscosity solution of (\ref{obstacleipde}) in the class $\pgc$.

\end{thm}
\no \proof : The proof is similar to the case without obstacle and based on the following facts:

(i) $u$ is continuous and belongs to $\pg$ ; 

(ii) The solution of the BSDE (\ref{reflec-bsde}) exists and is unique and is connected to $u$ by relation (\ref{fctu}); 

(iii) The characterization of the jumps of $Y^{t,x}$ by relation 
(\ref{caracsaut}) ; 

(iv) The generalization of Barles et al.'s definition (from the case without obstacle) to the case with obstacle and which is given in Appendix (Definition \ref{defipdeobst}). This generalization coincides with our definition when the generator $h$ does not depend on $\z$. 
\ms

\no The details of the proof are almost the same as the ones of the proof of Theorem \ref{mainthm} therefore they are left to the care of the reader. \qed

\section{Conclusion}
In this paper, we provide a new theoretical result of existence and uniqueness for solutions of some general class of non linear IPDEs. Especially and to our
 knowledge, there does not exist any study concerning viscosity solutions of such equations without assuming the monotonicity condition on the driver
 (with respect to its jump component).
 We note that our result deeply relies on the relationship between the solution of the non linear
 IPDE and the one of the related BSDE with jumps, relation given by the Feynman-Kac's formula.
 We also mention that since our proof is based on this relationship with some explicit BSDE (or reflected BSDE) with jumps,
 we obtain without additional difficulties the (existence and uniqueness) result both for the multidimensional case and for non linear IPDEs with one obstacle (see last Section \ref{finalremark}).
 As a consequence, this enlarges the class of economic and
 financial optimizations and/or control problems we can deal with which naturally 
 lead to the study of partial differential equations (or system of equations). Additionnaly
 , this study reinforces the interest of 
 using probabilistic tools in order to study PDEs or
 system of variational inequalities related to optimization problems.
 Another future application we have in mind is that this new result could be applied to obtain numerical implementations of such IPDEs.
 
\begin{center}
{\bf Appendix}
\end{center}  
  
\begin{defn}\label{bbpdef}\cite{BarlesBuckPardoux}: Barles et al.'s definition of a viscosity solution of (\ref{ipde-intro}). 
  \ms
    
\no We say that a family of deterministic functions $u =(u^i)_{i=1,m}$, defined on $\esp$ and $\r$-valued and such that for any $i\in \mo$, $u^i$ is continuous, is viscosity sub-solution (resp. super-solution)
 of the IPDE (\ref{ipde-intro}) if, for any $i\in \mo$: \\
(i) $\forall x\in \r^k$, $u_i(T, x) \le g_i(x) $ (resp. $u_i(T, x) \ge g_i(x) $) ;\\
 (ii)  For any  $(t,x)\in (0,T)\times \r^k$ and any function of class $\cC^{1,2}(\esp)$ such that $(t,x)$ is a global maximum point of 
$ u_i -\phi$ (resp. a global minimum point of $ u_i -\phi$) and $(u_i-\phi)(t,x)=0$, one has
$$ -\partial_t \phi(t,x)
-\mathcal{L}^X \phi(t,x) -h^{(i)}(t,x, (u^i(t,x))_{i=1,m}, \sigma^\top(t,x)D_x\phi(t,x), B_i(\phi)(t,x)) \le 0, 
$$
(resp. 
$$-\partial_t \phi(t,x)
-\mathcal{{L}}^X \phi(t,x) -h^{(i)}(t,x, (u^i(t,x))_{i=1,m}, \sigma^\top(t,x)D_x\phi(t,x), B_i(\phi)(t,x)) \ge 0.)
$$

The family $u =(u^i)_{i=1,m}$ is a viscosity solution of (\ref{ipde-intro}) if it is both a viscosity sub-solution and viscosity super-solution. 
\end{defn}
 
The adaptation of this definition to the case when there is an obstacle is the following (see 
\cite{HarrajOuknine05} or \cite{sulem}). 
  
  \begin{defn}\label{defipdeobst}
      
\no We say that a deterministic continuous functions $u$, defined on $\esp$ and $\r$-valued, is a viscosity sub-solution (resp. super-solution) of the IPDE (\ref{obstacleipde}) if: \\
(i) $\forall x\in \r^k$, $u_i(T, x) \le g(x) $ (resp. $u(T, x) \ge g(x) $) ;\\
 (ii)  For any  $(t,x)\in (0,T)\times \r^k$ and any function of class $\cC^{1,2}(\esp)$ such that $(t,x)$ is a global maximum point of 
$ u_i -\phi$ (resp. a global minimum point of $ u-\phi$) and $(u-\phi)(t,x)=0$, one has
$$\min\Big\{u(t,x)-\ell (t,x); -\partial_t \phi(t,x)
-\mathcal{L}^X \phi(t,x) -h(t,x, u(t,x), \sigma^\top(t,x)D_x\phi(t,x), B(\phi)(t,x))\Big\} \le 0, 
$$
(resp. 
$$\min\Big\{u(t,x)-\ell (t,x); -\partial_t \phi(t,x)
-\mathcal{L}^X \phi(t,x) -h(t,x, u(t,x), \sigma^\top(t,x)D_x\phi(t,x), B(\phi)(t,x))\Big\} \ge 0.)
$$

The function $u $ is a viscosity solution of (\ref{obstacleipde}) if it is both a viscosity sub-solution and viscosity super-solution. 
\end{defn}
\begin{rem} (i) If $h$ does not depend on $\zeta$ then Definitions \ref{od2} and \ref{defipdeobst} coincide ;
\ms

\no (ii) It is shown in 
\cite{HarrajOuknine05} or \cite{sulem} that when 
$h$ is non-decreasing $\wr$ $\z$ and $\g$ is moreover non-negative then the function $u$ defined in(\ref{fctu}) is the unique continuous viscosity solution of (\ref{obstacleipde}) in the sub-class of $\pg$ of continuous functions. \qed
\end{rem}

\end{document}